\documentclass[12pt]{article}  

\usepackage{amsmath,amsxtra,amssymb,latexsym, amscd}
\usepackage[mathscr]{eucal}

\begin{document}

\title{ \huge\bf Dedekind zeta function and BDS conjecture 
 }

\def\1{\rule{0cm}{0cm}} \def\qd{\rule{3mm}{3mm}} \def\BB{$\bullet$}
\def\bz{\boldsymbol{z}} 
\def\bx{\boldsymbol{x}} \def\by{\boldsymbol{y}} \def\bS{\boldsymbol{S}}
\def\I{{\cal I}} \def\A{{\cal A}} \def\C{\mathbb{C}}
\def\bv{\boldsymbol{v}}
\def\bn{\boldsymbol{n}} \def\by{\boldsymbol{y}} \def\bS{\boldsymbol{S}}
\def\H{{\mathcal H}} \def\U{{\cal U}} \def\D{{\cal D}}
\def\bc{\boldsymbol{c}} \def\bZh{\boldsymbol{Z}}
\def\ba{\boldsymbol{a}}\def\bu{\boldsymbol{u}}
\def\br{\boldsymbol{r}}\def\bzero{\boldsymbol{0}}
\def\R{\mathbb R}\def\K{\mathbb K}\def\Q{\mathbb Q}\def\F{\mathbb F}
\def\b{\mathfrak b} 
\def\a{\mathfrak a} \def\p{\mathfrak p} \def\f{\mathfrak f}
\def\C{\mathbb C}

\author{
{\bf Dang Vu Giang}\\
Hanoi Institute of Mathematics\\
Vietnam Academy of Science and Technology\\
18 Hoang Quoc Viet, 10307 Hanoi, Vietnam\\
       e-mail: dangvugiang@yahoo.com\\
\1\\
}
\maketitle 

\noindent {\bf Abstract.} 

\bigskip

\par\noindent{\sl Keywords:}     


\section{Introduction}

\noindent
Let $\K$ be an algebraic number field with degree $n=[\K:\Q]$.  Dedekind zeta function of $\K$ is first defined for complex numbers  $s$  with $\Re(s)>1$  by convergent  Dirichlet series
\[\zeta_\K (s) = \sum_{\a} |\a|^{-s}\]
where  $\a$ ranges through the non-zero ideals of the algebraic integers’ ring $\mathcal{O}_\K $  of $\K$ and  $ |\a| $ denotes the absolute norm of $\a$ which is  the cardinality of quotient ring $\mathcal{O}_\K/\a$. It is well known that $ |\a\b| =|\a| |\b|. $
for every ideals $\a$ and $\b$ of $\mathcal{O}_\K $.
Clearly, if $\K=\Q$,  the Dedekind zeta function $\zeta_\Q$ of $\Q$ becomes the classical Riemann zeta function $\zeta$. Moreover, for $\K=\Q(e^{2\pi i/m})$,
\[\zeta_{\Q(\exp\frac{2\pi i}m)}(s)=\zeta(s)\prod_{\chi}L(s,\chi)\prod_{\p|m}(1-|\p|^{-s})^{-1}.\]
 Here, $\chi$ denotes a non-principal Dirichlet character modulo $m$.
The Euler product of Dedekind zeta function of $\K$ is a product over all the prime ideals $\p$ of $\mathcal{O}_\K $  

\[
{1\over\zeta_\K (s)}= \prod_{\p} (1-|\p|^{-s})=\sum_{\a}{\mu_\K(\a)\over|\a|^s}
,\qquad\text{ for Re}(s)>1.
\]
This is the expression in analytic terms of uniqueness of  the prime factorization of an ideal  $\a$ in $\mathcal{O}_\K $ (Dedekind domain).  
The M\"obius function $\mu_\K$ in the field $\K$ is defined by
\[\mu_\K(\a)=\left\{\begin{matrix} 1 &\text{ if }& \a=\mathcal{O}_\K\\
(-1)^r&\text{ if }& \a=\p_1\p_2\cdots\p_r\\
0&\text{ if }& \p^2|\a\\
\end{matrix}
\right..
\] 
For example, $\mu_\K(5)=1$ for $\K=\Q(i)$, but $\mu_\Q(5)=\mu(5)=-1$.
For $\Re(s)>1$, $\zeta_\K (s) $ is given by a convergent product of infinite non-zero numbers,  so $\zeta_\K (s) $   is non-zero in this region. Erich Hecke first proved that $\zeta_\K (s) $  has an analytic continuation to the complex plane as a meromorphic function, having a simple pole only at $s=1$. Let $\Delta_\K$ denote the discriminant of $\K$,  $r_1$ and $r_2$   denote the number of real and complex places of $\K$, and 
\[
\Gamma_{\R}(s)=\pi^{-s/2}\Gamma(s/2),\quad 
\Gamma_{\C}(s)=2(2\pi)^{-s}\Gamma(s).\]
\[\Lambda_\K(s)=\left|\Delta_\K\right|^{s/2}
\Gamma_{\R}(s)^{r_1}
\Gamma_{\C}(s)^{r_2}\zeta_\K(s) \qquad
\Xi_\K(s)=\tfrac12(s^2+\tfrac14)\Lambda_K(\tfrac12+is).\]
We have the functional equation
\[\label{functional}
\Lambda_\K(s)=\Lambda_\K(1-s) \qquad
\Xi_\K (-s)= \Xi_\K (s).\]
Moreover,
\[\lim_{s\rightarrow1}(s-1)\zeta_\K(s)=\frac{2^{r_1+r_2} \pi^{r_2} h(\K)R(\K)}{w(\K)\sqrt{|\Delta_\K|}}.\]
In the virtue of this functional equation, it is also conjectured that the complex roots of $\zeta_\K (s) $   are on the critical line $\Re s=1/2$. This is called the extended Riemann hypothesis. If this hypothesis is true then every complex root of Riemann zeta function $\zeta$ and Dirichlet $L-$function $L(\chi)$ is on the line $\Re s=1/2$. 
The values of the Dedekind zeta function at integers encode important arithmetic data of the field $\K$. For example, the class number $h(\K)$ of $\K$  relates the residue at $s=1$, the regulator $R$($\K$) of $\K$, the number  $w$($\K$) of unity roots in $\K$, the  discriminant of $\K$, and the number of real and complex places of $\K$. From the functional equation and the fact that $\Gamma$  is infinite at all integers less than or equal to zero, it follows that  $\zeta_\K(s)$ vanishes at all negative even integers. It even vanishes at all negative odd integers unless  $\K$  is totally real number field. In the totally real case,  Carl Ludwig Siegel showed that $\zeta_\K(s)$  is a non-zero rational number at negative odd integers. Stephen Lichtenbaum conjectured specific values for these rational numbers in terms of the algebraic $K-$theory of  $\K$. Another example is at $s=0$ where $\zeta_\K(s)$   has a zero with order $r$ equal to the rank the unit group of $\mathcal{O}_\K $  and the leading term given by

\[\lim_{s\rightarrow0}s^{-r}\zeta_\K(s)=-\frac{h(\K)R(\K)}{w(\K)}.\]
Two fields are called arithmetically equivalent if they have the same Dedekind zeta function. There are examples  of non-isomorphic fields that are arithmetically equivalent. In particular some of these fields have different class numbers, so the Dedekind zeta function of a number field does not determine its class number. Now note that
\[\ln\zeta_\K(s)=\sum_{n=1}^\infty\sum_{\p}\frac1{n|\p|^{ns}},
\] 
\[\ln\zeta(s)+\ln\zeta(s-1)=\sum_{n=1}^\infty\sum_{p}\frac{p^n+1}{n|p|^{ns}}.
\]
For an elliptic curve $E$ over $\Q$ we define
\[\ln L(E,s)=\sum_{n=1}^\infty\sum_{p\not|\text{  } 2N}\frac{\#E(\F_{p^n})}{n|p|^{ns}}
\]
We have $\#E(\F_{p^n})=p^n+1$ if $E$ denotes the elliptic curves $y^2=x(x-N)(x+N)$ and $p=4k+3$.

{\footnotesize  

\bigskip

\par\noindent{\bf Acknowledgement.} Deepest appreciation is extended towards the NAFOSTED  (the National Foundation for Science and Techology Development in Vietnam) for the financial support.}


\begin{thebibliography}{999}
\bibitem{Gasper 2} G. Gasper, Using integrals of squares to prove
that certain real-valued special functions to prove that the P\'{o}lya
$\Xi^{*}\left(z\right)$ function, the functions $K_{iz}\left(a\right),\ a>0$
and some other entire functions have only real zeros, dedicated to
Dan Waterman on the occasion of his 80th birthday.
\end{thebibliography}
\end{document}